\documentclass[12pt,a4paper]{amsart}

\makeatletter
\long\def\unmarkedfootnote#1{{\long\def\@makefntext##1{##1}\footnotetext{#1}}}
\makeatother

\usepackage{amsmath,amssymb}%,a4}
\usepackage{enumerate,amssymb}
\usepackage{color}
\usepackage{tikz}
\usetikzlibrary{arrows}
\theoremstyle{plain}
\newtheorem{thm}{Theorem}[section]

\newtoks\prt
\newtheorem{proclaim}[thm]{\the\prt}
\theoremstyle{definition}

\def\eqn#1$$#2$${\begin{equation}\label#1#2\end{equation}}

\numberwithin{equation}{section}

\def\diam{\operatorname{diam}}
\def\dist{\operatorname{dist}}

\def\epsilon{\varepsilon}
\def\en{\mathbb N}
\def\er{\mathbb R}

\def\J{\mathcal{J}}

\def\mira{\mathcal L_n}
\def\mir1{\mathcal L_1}

\def\Q{\mathcal{Q}}
\def\r{\widetilde{r}}

\def\phi{\varphi}
\def\sgn{\operatorname{sgn}}

\def\rn{\mathbb R^n}
\def\spt{\operatorname{supp}}
\def\sgn{\operatorname{sgn}}

\def\V{\mathbb V}

\def\z{\widetilde{z}}

\def\ve{\boldsymbol v}
\def\w{\boldsymbol w}

\newtoks\by
\newtoks\paper
\newtoks\book
\newtoks\jour
\newtoks\yr
\newtoks\pages
\newtoks\vol
\newtoks\publ
\def\ota{{\hbox\vol{???}}}
\def\cLear{\by=\ota\paper=\ota\book=\ota\jour=\ota\yr=\ota
\pages=\ota\vol=\ota\publ=\ota}
\def\endpaper{\the\by, {\the\paper},
\textit{\the\jour} \textbf{\the\vol} (\the\yr), \the\pages.\cLear}
\def\endbook{\the\by, \textit{\the\book}, \the\publ.\cLear}
\def\endprep{\the\by, \textit{\the\paper}, \the\jour.\cLear}
\def\endyearprep{\the\by, \textit{\the\paper}, \the\jour, (\the\yr).\cLear}
\def\name#1#2{#2 #1}
\def\nom{ \rm no. }

\title{Ciarlet Ne\v{c}as condition in fractional Sobolev spaces}

\author{$^1$Stanislav Hencl}

\author{$^2$Jarom\'ir Mielec}

\author{$^3$Kaushik Mohanta}

\address{$^{1,2}$Department of Mathematical Analysis, Charles University,
	So\-ko\-lovsk\'a 83, 186~00 Prague 8, Czech Republic}
	
\address{$^3$School of Mathematics and Computer Science, Indian Institute of Technology Goa, Ponda-403401, India}

\email{\tt $^1$hencl@karlin.mff.cuni.cz,  $^2$mielec@karlin.mff.cuni.cz, $^3$kaushik@iitgoa.ac.in}

\date{\today}

\thanks{The authors were supported by the grant GA\v{C}R P201/24-10505S}

\begin{document}

\begin{abstract}
Let $s\in(\frac{n}{n+1},1)$, $\Omega\subset\rn$ be an open set and let $f\in W^{s,n/s}(\Omega,\rn)$ be mapping with positive distributional Jacobian $\J_f>0$ which models some deformation in fractional Nonlinear Elasticity. We show change of variables formula in this class and as a consequence we show that the analogue of Ciarlet-Ne\v{c}as condition $\J_f(\Omega)=|f(\Omega)|$ implies that our mapping is one-to-one a.e.   
\end{abstract}

\maketitle

\section{Introduction}

Let $\Omega\subset\rn$ be an open set and let $f\colon\Omega\to\rn$ be a mapping.
In this paper, we study classes of mappings $f$ that might serve as deformations in some models in Nonlinear Elasticity. Following the pioneering papers of Ball \cite{Ball} and Ciarlet and Ne\v{c}as \cite{CN} we ask if our mapping is in some sense injective as the physical `non-interpenetration of the matter' asks a deformation to be 
one-to-one. One of the standard possibilities is to use the Ciarlet-Ne\v{c}as condition which works for $W^{1,p}$, $p>n$, spaces and asks mappings with $J_f\geq 0$ to satisfy 
$$
\int_{\Omega} J_f= |f(\Omega)|.
$$
Mappings in this class are one-to-one a.e. which is not difficult to show using the standard change of variables formula $\int_{\Omega} |J_f|=\int_{\rn}N(f,y,\Omega)\; dy$. Moreover, this condition is stable under weak convergence and hence can be used easily in the direct method of the Calculus of Variations. 
The main aim of this paper is to introduce some variant of the Ciarlet-Ne\v{c}as condition in the context of fractional Nonlinear Elasticity. 
As far as we know this is a first result in this area which allows us to show injectivity a.e. of the deformations. 

Recently there has been a growing interest in nonlocal models of Nonlinear Elasticity, especially since the introduction of the peridynamics model by Silling \cite{Si}. In those models we do not work with the derivative but with fractional derivatives instead (see e.g. \cite{BCM}, \cite{KS} and \cite{BCM2} and references given there for some recent results) as it models well some nonlocal interactions and allow for singularities. 
In our approach we need to work with some variant of the Jacobian which allows us to show some variant of the change of variables formula. There is no reasonable pointwise variant of the derivative which could be used in fractional Sobolev spaces for change of variables but we can use the idea of distributional Jacobian introduced by Brezis and Nguyen \cite{BN}. They showed that in fractional order Sobolev spaces $W^{\frac{n-1}{n},n}$ we can define the distributional Jacobian as
$$
\J_f(\phi)=\lim_{k\to\infty}\int_{\Omega} J_{f_k}(x)\phi(x)\; dx\text{ for every }\phi\in C_c^1(\Omega),  
$$
where $f_k$ are smooth approximations of $f$ in $W^{\frac{n-1}{n},n}$.
The existence of this limit was extended in Li and Schikorra in \cite{LS} to every $s>\frac{n}{n-1}$, $f\in W^{s,n/s}$ and  $\phi \in W^{(1-s)n,\frac{1}{1-s}}$. Moreover, it was shown in \cite{LS} that mappings in this class for $s>\frac{n}{n+1}$ with $\J_f>0$ are continuous and sense preserving and these properties are crucial and natural when we model deformations. Note that $\J_f$ can be represented by a measure in this case as it is a nonnegative distribution. 
%We would like to study further properties of these maps like validity of change of variables formula or injectivity a.e. 

In our main result we first show some variant of change of variables formula for the distributional Jacobian and using it we show that each mapping with $\J_f>0$ which satisfies the Ciarlet-Ne\v{c}as condition is one-to-one a.e.

%Our main result is the following (WE HOPE). EXPLAIN THE DEFINITION OF $f(\partial \Omega)$ 

\prt{Theorem}
\begin{proclaim}\label{main}
Let $n\geq 2$, $s\in(\frac{n}{n+1},1)$ and $\Omega\subset\rn$ be open with smooth boundary. 
Assume that $f\in W^{s,n/s}(\Omega,\rn)$ satisfies $\J_f>0$. 
Then we have 
\eqn{change}
$$
\J_f(\psi\circ f)=\int_{\rn}\deg(f,\Omega,y)\; \psi(y)\; dy
\text{ for every }\psi\in C^1_c(\rn\setminus f(\partial \Omega)).
$$

Let $U\subset\subset\Omega$ be a domain such that $|f(\partial U)|=0$, $|f^{-1}(f(\partial U))|=0$ and that $f$ satisfies the Ciarlet Ne\v{c}as condition 
\eqn{CN}
$$
\J_f(U)=\bigl|f(U)\bigr|
$$
where the symbol $\J_f(U)$ denotes the measure of the set $U$. 
Then for a.e. $y\in f(U)$ there is only one $x\in U$ with $f(x)=y$. 
%Then $f$ has $\deg(f,\Omega,y)=1$ for every $y\in f(\Omega)\setminus f(\partial \Omega)$. 
%It follows that for every $y\in f(\Omega)\setminus f(\partial \Omega)$ there is either only one $x\in\Omega$ with $f(x)=y$ or there are uncountably many such $x$. 
\end{proclaim}

Let us note that once we prove change of variables formula \eqref{change} it is not difficult to show that Ciarlet-Ne\v{c}as condition implies that our mapping satisfies that its topological degree equals to 1 a.e. However this does not imply that the mapping is 1-1 a.e. as stated above. Indeed, in \cite{BHM} you can find an example of $W^{1,n-1}$ mapping (strong limit of homeomorphisms) where the degree is 1 a.e., but the preimage of points is a continuum for a set of positive measure. We show using some "capacitary type" estimates that this cannot happen for our $W^{s,n/s}$ mapping and our map is in fact even 1-1 a.e.

In the statement of the Ciarlet-Ne\v{c}as condition \eqref{CN} we assume that we are on some subdomain so that the continuous representative of $f$ (given by the result \cite{LS}) is well defined on $\partial U$. Alternatively it would be possible to assume that $U$ and $f(U)$ have smooth boundary and to extend $f$ over the boundary in a standard way by some mirroring (see e.g. the approach in \cite{M}).

%The last conclusion with uncountably many preimages can really happen. We can easily construct a Lipschitz mapping which squeezes a segment to a point and it is a homeomorphism elsewhere. Moreover, there exists a continuous sense-preserving mapping in $W^{1,n-1}$ such that for every $y$ in a Cantor type set of positive measure we know that $f^{-1}(y)$ is a continuum (see \cite[Theorem 1.1]{BHM}). It is not clear if such a mapping could be created with $W^{s,\frac{n}{s}}$ regularity. 

In the usual change of variables formula for smooth mappings we have $\int_{\Omega}J_f(x) \psi(f(x))\; dx$ on the left hand side of \eqref{change}. 
This means that the measure $\J_f$ has no singular part and we have the validity of the Lusin $(N)$ condition, i.e. sets of measure zero are mapped to sets of measure zero. This key property may unfortunately fail in our situation in any $W^{s,p}$ for any $s<1$. 

\prt{Theorem}
\begin{proclaim}\label{exampleJarda}
Let $n\geq 2$, $s\in(0,1)$ and $p\in(1,\infty)$. Then there is a homeomorphism $f\in W^{s,p}([-1,1]^n,\rn)$ such that $f(x)=x$ for $x\in\partial [-1,1]^n$ which fails the Lusin $(N)$ condition, i.e. there is $N\subset (-1,1)^n$ such that $|N|=0$ but $|f(N)|>0$. 
\end{proclaim}

Let us note that the assumption about the positivity (or nonnegativity) of the distributional Jacobian is really needed in this theory as the pointwise Jacobian does not need to carry any information about 
the topological behavior of the mapping. This is shown e.g. by the following example which can be constructed in any $W^{s,p}$ for any $s<1$. 

\prt{Theorem}
\begin{proclaim}\label{exampleKaushik}
Let $n\geq 2$, $s\in(0,1)$ and $p\in(1,\infty)$. Then there is a continuous mapping $\tilde{f}\in W^{s,p}([-1,1]^n,\rn)$ such that classical derivative $\nabla \tilde{f}$ exists a.e. and $J_{\tilde{f}}(x)=\det \nabla \tilde{f}(x)>0$ a.e. but $\tilde{f}$ is not sense preserving as $\tilde{f}(x)=[-x_1,x_2,\hdots, x_n]$ for $x\in\partial [-1,1]^n$. 
\end{proclaim}

%FURTHER QUESTIONS - IF $sp>0$ AND $\J_F\geq 0$, DOES IT FOLLOW THAT F IS CONTINUOUS AND SENSE PRESERVING? THEN WE COULD DO CALC VAR THERE AS THE CLASS IS WEAKLY CLOSED. Note that it is enough to have $\J_f\geq 0$ in \cite[Prop 1.7 and 1.8]{LS}. Do we really need it in \cite[Prop 1.9]{LS}??? Proof of \cite[Prop 1.7 and 1.8]{LS} for continuous $f$ is actually not difficult: Given $y\in f(\Omega)\setminus f(\partial\Omega)$ find nonnegative smooth $\psi\approx \delta_y$. As $\deg$ is constant on the support of $\psi$ we obtain by \eqref{change} that 
%$$
%\deg(f,\Omega,y)=\J_f(\psi\circ f).
%$$
%Since $\psi\geq 0$ we  easily get $\deg\geq 0$ actually for any subdomain. 

\section{Preliminaries}

By $Q(c,r)$ we denote an open cube centered at $c\in\rn$ with sidelength $2r$. %For $\alpha>0$ and a cube $Q=Q(c,r)$ we denote $\alpha Q=Q(c,\alpha r)$. 

By $C_c^1(\Omega)$ we denote the class of function with continuous first order derivatives and with compact support in $\Omega$.

We need the following lemma about differentiability of mappings in special form from \cite[Lemma 2.1]{HK}. 

\prt{Lemma}
\begin{proclaim}\label{radial}
Let $\rho:(0,\infty)\to(0,\infty)$ be a strictly monotone function and $\rho\in C^1((0,\infty))$. Then for the mapping
$$f(x)=\frac{x}{|x|}\rho(|x|), \quad x\neq 0$$
we have for almost every $x$
$$|Df(x)|= \max\Bigl\{\frac{\rho(|x|)}{|x|}, |\rho'(|x|)| \Bigr\}\text{ and }
 J_f(x)= \rho'(|x|) \Bigl(\frac{\rho(|x|)}{|x|}\Bigr)^{n-1}\ .$$
 %$$\text{ and }|\adj Df(x)|\sim \max\Bigl\{\frac{\rho(\|x\|)}{\|x\|}, |\rho'(\|x\|)| \Bigr\} 
 %\Bigl(\frac{\rho(\|x\|)}{\|x\|}\Bigr)^{n-2}.$$
\end{proclaim}

\subsection{H\"older spaces}
The space of $\alpha$-H\"older continuous functions over a bounded domain $\Omega$ is a Banach space when equipped with the norm 
$$
\|\cdot\|_{C^{0,\alpha}(\Omega)}:= \|\cdot\|_{L^\infty(\Omega)} +|\cdot|_{C^{0,\alpha}(\Omega)},
$$
where
$$
|f|_{C^{0,\alpha}(\Omega)}:=\sup_{x,y\in \Omega,\ x\neq y} \frac{|f(x)-f(y)|}{|x-y|^{\alpha}}.
$$
We mention the following result from \cite[Lemma~6.33]{GT}.
\prt{Lemma}
\begin{proclaim}\label{cpt-embedding-holder}
Let $\Omega$ be a bounded open set in $\rn$ and $0<\beta<\alpha\leq 1$. Then $C^{0,\alpha}(\Omega)$ is compactly embedded inside $C^{0,\beta}(\Omega)$.
\end{proclaim}

We know from (\cite[Theorem 8.2]{DiPP}) that for $sp>n$ we get the embedding into H\"older spaces with $\alpha=s-\frac{n}{p}$
%we get embedding into H\"older spaces on this $(n-1)$-dimensional square $Q$ (\cite[Theorem 8.2]{DiPP}) with a constant $C$ independent on the size of the cube $Q$.
$$
\sup_{x,y\in Q(0,1)} \frac{|f(x)-f(y)|}{|x-y|^{\alpha}}\leq C\Bigl(\|f\|_{L^p(Q(0,1))}+\int_{Q(0,1)}\int_{Q(0,1)}\frac{|f(x)-f(y)|^p}{|x-y|^{n+sp}}\; dx\; dy\Bigr). 
$$
By a simple scaling argument (i.e. change of variables $y=Rx$) we get
$$
R^{\alpha} \sup_{x,y\in Q(0,R)} \frac{|f(x)-f(y)|}{|x-y|^{\alpha}}\leq C\Bigl(R^{-n/p}\|f\|_{L^p(Q(0,R))}+R^{\alpha} \int_{Q(0,R)}\int_{Q(0,R)}\frac{|f(x)-f(y)|^p}{|x-y|^{n+sp}}\; dx\; dy\Bigr). 
$$
By limiting $R\to\infty$ we get that the term $\|f\|_{L^p(Q(0,R))}$ is redundant - more precisely we fix $f\in W^{s,p}(Q(0,R_0))$, extend it to $f\in W^{s,p}(Q(0,2R_0))$ with compact support and we use the previous inequality on this $f$ for $R\geq R_0$ and by sending $R\to\infty$ we see after little work that the $L^p$ term is redundant. 
It follows that we have the embedding
\eqn{emb}
$$
\sup_{x,y\in Q(0,R)} \frac{|f(x)-f(y)|}{|x-y|^{\alpha}}\leq C \Bigl(\int_{Q(0,R)}\int_{Q(0,R)}\frac{|f(x)-f(y)|^p}{|x-y|^{n+sp}}\; dx\; dy\Bigr)^{\frac{1}{p}}
$$
and the scaling tells us that the constant on the right-hand side does not depend on $R$.

%|f(x_Q)-f(y_Q)|\leq C |x_Q-y_Q|^{s-\frac{n-1}{p}}\Bigl(\int_{Q}|f|^p+\int_{Q}\int_Q\frac{|f(x)-f(y)|^p}{|x-y|^{n-1+sp}}\; dx\; dy\Bigr)^{\frac{1}{p}}.
%$$

\subsection{Distributional Determinant}

The following statements were shown in \cite{LS} (see the paper for exact definitions):

\prt{Theorem}
\begin{proclaim}\label{defjac}
Let $n\geq 2$, $s\in(\frac{n-1}{n},1)$ and $\Omega\subset\rn$ be open with smooth boundary. Assume that $f\in W^{s,n/s}(\Omega,\rn)$ and that $\phi\in W^{(1-s)n,\frac{1}{1-s}}(\Omega)$. 
Then the following extension of the Jacobian operator
$$
\J_f(\phi)=\lim_{k\to\infty}\int_{\Omega}J_{f_k}(x)\phi_k(x)\; dx
$$
is well-defined for any $f_k\in C^{\infty}(\overline{\Omega})$ which is a smooth approximation of $f\in W^{s,n/s}(\Omega,\rn)$ and any $\phi_k\in C^{\infty}(\overline{\Omega})$ which is a smooth approximation of $\phi\in W^{(1-s)n,\frac{1}{1-s}}(\Omega)$. 
\end{proclaim}

\prt{Theorem}
\begin{proclaim}\label{sense}
Let $n\geq 2$, $s\in(\frac{n}{n+1},1)$ and $\Omega\subset\rn$ be open with smooth boundary. Assume that $f\in W^{s,n/s}(\Omega,\rn)$ satisfies $\J_f>0$. Then $f$ is sense-preserving, weakly pseudomonotone and continuous. 
\end{proclaim}

The main idea of the proof of previous theorem in \cite{LS} is to try to obtain the "analogy" of the definition of the degree as
$$
\deg(f, B, y)=\int_{\Omega}J_f(x)\psi(f(x))\; dx
$$
where $\psi$ is the smooth approximation of Dirac measure $\delta_y$. In order to do so we need to have finiteness of the above integral, i.e. we need to test the distributional Jacobian by the $f$ itself (more precisely by some function that has $f$ inside). This is exactly what forces the restriction $s>\frac{n}{n+1}$ as they can define $\J_f(\phi)$ for every $f\in W^{s,n/s}$ and $\phi\in W^{(1-s)n,1/(1-s)}$. 

The following result is taken from \cite[Lemma A.2]{LS}. We have modified it slightly to suit our purpose; still, the same proof as in \cite{LS} works for this version.
\prt{Lemma}
\begin{proclaim}\label{lm-2}
	Let $\Omega\subset \rn$ be a bounded domain with smooth boundary, $s\in (0,1)$, $p>1$, and $f\in W^{s,p}(\Omega)$. Denote, $H_r:=\{x\in \rn\ |\ \pi_i(x)=r\}$ where $\pi_i$ is the projection to $x_i$ axis. Then $f\in W^{s,p}(H_r\cap \Omega)$ for a.e. $r$ and we have 
	$$
	\int_{\mathbb{R}} \|f\|_{W^{s,p}(H_r\cap \Omega)}^p dr
	\leq C(\Omega, n,p, s) \|f\|_{W^{s,p}(\Omega)}^p.
	$$
\end{proclaim}

The following result is taken from \cite[Theorem~3]{BN}.

\prt{Theorem}
\begin{proclaim}\label{BN-theorem}
Let $n \geq 2$. For all $f, \ g \in W^{\frac{n-1}{n}, n}(\rn)$, and for all $\psi \in
C_c^1(\rn)$, we have
$$
|\J_f(\psi) - \J_g(\psi)|
\leq  C(n,\Omega) \|f - g\|_{W^{\frac{n-1}{n},n}} \left( \|f\|_{W^{\frac{n-1}{n},n}}^{n-1}(\Omega) +\|g\|_{W^{\frac{n-1}{n},n}}^{n-1}(\Omega)\right) \|\nabla \psi\|_{L^\infty(\rn)}.
$$
\end{proclaim}

\subsection{Degree}

Given a continuous mapping $f$ from $\Omega\subset \rn$ to $\rn$ and $y\in\rn\setminus f(\partial\Omega)$, the topological degree $\deg(f,\Omega,y)$ describes the number of preimages $f^{-1}(y)$ in $\Omega$ taking the orientation into consideration. For a smooth map, we can define the topological degree as
\[\deg (f,\Omega, y)= \sum_{\{x\in\Omega:f(x)=y\}}\sgn (J_f(x))\]
if $J_f(x)\neq 0$ for each $x\in f^{-1}(y)$. This definition can be extended to arbitrary continuous mapping and each point in $\rn\setminus f(\partial\Omega)$ (see e.g. \cite{FG}).

A continuous mapping $f:\Omega\subset\rn\to\rn$ is called \emph{sense-preserving}, if 
\[\deg(f,\Omega^\prime, y)>0\]
for all domains $\Omega^\prime\subset\subset \Omega$ and all $y\in f(\Omega^\prime)\setminus f(\partial\Omega^\prime)$.

%We need the following observation about continuous sense-preserving mappings - STANDA BELIEVES THAT IT IS TRUE BUT WE NEED TO FIND THE REFERENCE

%\prt{Lemma}
%\begin{proclaim}\label{imagepositive}
%    Let  $\Omega\subset\rn$ be and open set, let $f\in C(\Omega,\rn)$ be sense preserving and let %$U\subset\subset\Omega$ be a domain. Then $f(U)$ is exactly the bounded set of $\rn$ surrounded by %$f(\partial U)$, more precisely it contains the whole interior of this set and it may contain some %part of $f(\partial U)$. 
%\end{proclaim}

The following statements are also from the book \cite{FG}. 

\prt{Theorem}
\begin{proclaim}\label{decomp}
    Let  $f\in C(U,\rn)$ and let $y\notin f(\partial U).$ If $U_1,U_2\subset U$ are open and disjoint and $y\notin f(U\setminus (U_1\cup U_2))$ then
    \[\deg(f,U,y)=\deg(f,U_1,y)+\deg(f,U_2,y).\]
\end{proclaim}

\prt{Theorem}
\begin{proclaim}\label{aprox}
Let $f\in C(\Omega,\rn)$ and $y\notin f(\partial\Omega)$. Then for every $g\in C(\Omega,\rn)$ such that $\|f-g\|_{\infty}<\dist(y, f(\partial\Omega))$ we have
\[\deg(f,\Omega, y)=\deg(g,\Omega, y).\]
\end{proclaim}

We also have a change of variables via the degree.

\prt{Theorem}
\begin{proclaim}\label{var}
Let $n\geq 2$ and $\Omega\subset\rn$ be an open set. Assume that $f\in C^1(\Omega,\rn)$. Further assume that $G\subset\Omega$ is bounded and $\psi:\rn\to \er$ is a Borel measurable bounded function. Then 
$$
\int_{G} J_f(x)\psi(f(x))\; dx=\int_{\rn}\deg(f,G,y)\psi(y)\; dy.
$$
%and for $\psi\geq 0$ (???) we have
%$$
%\int_{\Omega} |J_f(x)|\psi(f(x))\; dx=\int_{\rn}N(f,\Omega,y)\psi(y)\; dy.
%$$
\end{proclaim}

\section{Change of Variables and Ciarlet Ne\v{c}as condition}

\begin{proof}[Proof of Theorem \ref{main}] 
	
	Recall from Theorem \ref{sense} that our mapping is continuous and sense-preserving. 

We take $f_k$ e.g. as the usual convolution approximation of $f$ with $f_k\to f$ in $W^{s,\frac{n}{s}}$ - note that $\Omega$ has smooth boundary and hence we can extend $f$ outside of $\Omega$ in the same class and $f_k$ are well-defined on $\overline{\Omega}$. Since $s>\frac{n}{n+1}$ we obtain from \cite[Theorem~2.1]{Jawe} that
$$
W^{s,\frac{n}{s}}\hookrightarrow W^{(1-s)n,1/(1-s)}\text{ for }s>\frac{n}{n+1}.
$$  
This implies that $f_k\to f$ in $W^{(1-s)n,1/(1-s)}$ and by the composition $\phi_k=\psi\circ f_k$ to $\phi=\psi\circ f$ in $W^{(1-s)n,1/(1-s)}$ as $\psi\in C_c^1$ (see \cite[Theorem~7]{BouSic}). 
By Theorem \ref{defjac} we obtain that
\eqn{a}
$$
\J_f(\psi\circ f)=\lim_{k\to\infty}\int_{\Omega}J_{f_{k}}(x) \psi(f_k(x))\; dx. 
$$

Now by the change of variables Theorem \ref{var} we obtain
\eqn{b}
$$
\int_{\Omega}J_{f_k}(x)\psi(f_k(x))\; dx=\int_{\rn}\deg(f_k,\Omega,y)\psi(y)\; dy. 
$$
Finally $\psi\in C^1_c(\rn\setminus f(\partial \Omega))$ so $\dist(\spt{\psi},f(\partial \Omega))>0$. Since $f_k$ are converging uniformly, we can find $k_0$ such that for all $k\geq k_0$ we have $\|f-f_k\|_{\infty}<\dist(\spt{\psi},f(\partial \Omega))$. It follows from Theorem \ref{aprox} that

$$
\deg(f_k,\Omega,y)=\deg(f,\Omega,y)\text{ for every }y\in \spt{\psi}\text{ and every }k\geq k_0. 
$$
This, together with \eqref{a} and \eqref{b}, implies \eqref{change}. 

Now it is not difficult to show that $\deg=1$ a.e. in the target. 
Analogously to \eqref{change} we obtain that
$$
\J^U_f(\psi\circ f)=\int_{\rn}\deg(f,U,y)\; \psi(y)\; dy
\text{ for every }\psi\in C^1_c(\rn\setminus f(\partial U)),
$$
where $\J^U_f$ is the restriction of $\J_f$ to functions supported in $U$. 
Using $|f^{-1}(f(\partial U))|=0$ we take $\psi_k\nearrow \chi_{f(U)}$ with $\psi_k\circ f|_U\nearrow \chi_{U}$ and passing to a limit in previous equality we obtain
$$
\J_f(U)=\J^U_f(U)=\int_{f(U)}\deg(f,U,y)\; dy
$$
where the first equality follows from the fact that $\J_f$ is a nonnegative measure. 
Now $|f(\partial U)|=0$ and Ciarlet Ne\v{c}as condition imply
$$
\int_{f(U)}\deg(f,U,y)\; dy=\int_{\rn}\deg(f,U,y)\; dy=\J_f(U)=|f(U)|
$$
which together with the fact that $f$ is sense preserving (thus $\deg\geq 1$ on $f(U)\setminus f(\partial U)$) easily imply that $\deg(f,U,y)=1$ for almost every $y\in f(U)\setminus f(\partial U)$. Since $\deg$ is locally constant we obtain that  $\deg(f,U,y)=1$ for every $y\in f(U)\setminus f(\partial U)$. 

Since $\deg(f,U,y)=1$ for such $y\in f(U)\setminus f(\partial U)$ we obtain that each such a point has at least one preimage. Assume that $f^{-1}(y)$ contains two points $a,b\in U$. We claim that $f^{-1}(y)$ intersects any hyperplane $H$ between $a$ and $b$ and thus $f^{-1}(y)$ is uncountable. Assume for contrary that 
$y\cap f(H)=\emptyset$. Then $H$ divides $U$ into two open sets $U_a\ni a$ and 
$U_b\ni b$ and $y\notin f(\partial U_a)\cup f(\partial U_b)\subset f(\partial U)\cup f(H)$. Since $f$ is sense preserving we obtain that $\deg(f,U_a,y)\geq 1$ and $\deg(f,U_b,y)\geq 1$. However this contradicts the additivity of the degree Theorem \ref{decomp}
\eqn{spor}
$$
1=\deg(f,U,y)=\deg(f,U_a,y)+\deg(f,U_b,y)\geq 2
$$
as $y\notin f(U\setminus (U_a\cup U_b))=f(H)$. 

Finally, we need to show that the number of preimages could be infinite $\# f^{-1}(y)=\infty$ only on a set of zero measure.
Assume for contrary that this is not the case. Then 
$$
Y:=\{y:\ \# f^{-1}(y)=\infty\}\text{ satisfies }|Y|>0.
$$

For $u\in\mathbb{S}^{n-1}$ and $t\in\er$ we denote by $H_u^t$ an $n-1$ dimensional hyperplane with a unit normal vector $u$
\[H_u^t=\Bigl\{x\in\rn:\sum_{i=1}^nu_ix_i=t\Bigr\}.\]
We define a set
\[A=\left\{(y,u,t)\in Y\times\mathbb{S}^{n-1}\times\er:f^{-1}(y)\cap H_u^t\neq\emptyset\right\}.\]
Note that $A$ is a measurable set as for each $K\subset Y$ compact the sets 
$$
A_K=\left\{(y,u,t)\in A:y\in K\right\}
$$ 
are closed sets. We fix $y\in Y$. It has (more than) two preimages $a\in\Omega$ and $b\in\Omega$, $a\neq b$. Now for each $u\in\mathbb{S}^{n-1}$ such that $u$ and $b-a$ are not orthogonal vectors there exists $t_a$ such that $a\in H_u^{t_a}$ and $t_b\neq t_a$ such that $b\in H_u^{t_b}$. Analogously to the paragraph before \eqref{spor} we can show that for each $t\in(t_a,t_b)$ we have $(y,u,t)\in A$ and hence 
$$
\int_\er \chi_A(y,u,t)dt>0.
$$ 
%Then $H_u^t$ divides $\Omega$ into two open sets $\Omega_a\ni a$ and $\Omega_b\ni b$ for each $t\in(t_a,t_b)$. It follows that such a triplet $(y,u,t)\in A$. 
Since for each $y\in Y$ there is a full measure of $u\in\mathbb{S}^{n-1}$ that are not orthogonal to $b-a$ we have
\[
\int_{\rn}\left(\int_{\mathbb{S}^{n-1}}\left(\int_\er \chi_A(y,u,t)dt\right)du\right)dy>0.
\]

By Fubini theorem we also have that 
\[
\int_{\mathbb{S}^{n-1}}\left(\int_\er\left(\int_{\rn} \chi_A(y,u,t)dy\right)dt\right)du>0.
\]
and thus there exist $u\in\mathbb{S}^{n-1}$ %such that 
%\[
%\int_\er\left(\int_{\rn} \chi_A(y,u,t)dy\right)dt>0.
%\]
and $T\subset\er$, $|T|>0$ such that
\eqn{aaa}
$$
\int_{\rn}\chi_A(y,u,t)\; dy>0\text{ for each }t\in T.
$$
Using Lemma \ref{lm-2} we can moreover fix $t\in T$ such that $f\in W^{s,p}(H_u^t\cap U)$. 
For this fixed $u$ and $t$ we have \eqref{aaa} and thus
$$
\bigl|f\bigl(f^{-1}(Y)\cap H_u^t\bigr)\bigr|>0.
$$

Assume without loss of generality that $H_u^t\cap U$ is a unit square which we divide into $N^{n-1}$ squares of sidelength $1/N$. We denote the collection of these squares $\Q$. Then we have
\[
0<\bigl|f\bigl(f^{-1}(Y)\cap H_u^t\bigr)\bigr|\leq \Bigl|\bigcup_{Q\in\Q}f(Q)\Bigr|\leq\sum_{Q\in\Q}|f(Q)|\leq \sum_{Q\in\Q}(\mathrm{osc}_Qf)^n.
\]
For each such a $(n-1)$-dimensional square $Q\in\Q$ we pick points $x_Q,y_Q\in Q$ where the oscillation of $f$ is  almost attained and we have the following
\begin{equation}\label{eee}
    0<\frac{\bigl|f\bigl(f^{-1}(Y)\cap H_u^t\bigr)\bigr|}{2}\leq  \sum_{Q\in\Q}|f(x_Q)-f(y_Q)|^n.
\end{equation}

Since $sp=s\frac{n}{s}=n>n-1$ we get embedding into H\"older spaces on this $(n-1)$-dimensional square $Q$ with a constant $C$ independent on the size of the cube $Q$ by \eqref{emb}
$$
|f(x_Q)-f(y_Q)|\leq C |x_Q-y_Q|^{s-\frac{n-1}{p}}\Bigl(\int_{Q}\int_Q\frac{|f(x)-f(y)|^p}{|x-y|^{n-1+sp}}\; dx\; dy\Bigr)^{\frac{1}{p}}.
$$
We raise this estimate to $n$, use $|x_Q-y_Q|\leq C/N$, apply H\"older's inequality and sum over $Q$ to get
$$
\begin{aligned}
\sum_{Q\in\Q} |f(x_Q)-f(y_Q)|^n &\leq \frac{C}{N^{(s-\frac{n-1}{p})n}} \sum_{Q\in\Q} \Bigl(\int_{Q}\int_Q\frac{|f(x)-f(y)|^p}{|x-y|^{n-1+sp}}\; dx\; dy\Bigr)^{s}\\
&\leq \frac{C}{N^{(s-\frac{n-1}{n}s)n}} \Bigl(\sum_{Q\in\Q} \int_{Q}\int_Q\frac{|f(x)-f(y)|^p}{|x-y|^{n-1+sp}}\; dx\; dy\Bigr)^{s}\Bigl(\sum_{Q\in\Q} 1\Bigr)^{1-s}\\
&\leq \frac{C}{N^s} \Bigl(\int_{H_u^t\cap\Omega}\int_{H_u^t\cap\Omega}\frac{|f(x)-f(y)|^p}{|x-y|^{n-1+sp}}\; dx\; dy\Bigr)^{s}\bigl(N^{n-1}\bigr)^{1-s}.\\
\end{aligned}
$$  
Since $(n-1)(1-s)<s$ (by our assumption $s>\frac{n}{n+1}$) we have that the right hand side tends to zero as $N\to \infty$ (recall that $f\in W^{s,p}(H_u^t\cap\Omega)$) which gives us a contradiction with (\ref{eee}).
\end{proof}

\subsection{Stability of Ciarlet Ne\v{c}as condition under weak convergence}

\prt{Lemma}
\begin{proclaim}\label{lm-1}
Let $\Omega\subset \rn$ be a bounded domain, let $U\subset\subset\Omega$ be a domain and let $s\in \left(\frac{n-1}{n},1\right)$. Assume that a sequence of mappings $\{f_k\}_k$, converges to $f$ weakly in $W^{s,\frac{n}{s}}(\Omega,\rn)$, $|f(\partial U)|=0$,  and assume that $f$ and that $ f_k$ are continuous and sense-preserving for all $k$. Then
\eqn{haha}
$$
\lim_{k\to\infty}|f_k(U)|\to |f(U)|.
$$
\end{proclaim}
\begin{proof}
%By Theorem \ref{sense} we know that $f$ and all $f_k$ are continuous and sense-preserving. 	
	
%We know by Lemma~\ref{lm-2} that
%$$
%\int_{\mathbb{R}} [f]_{W^{s,p}(H_r\cap \Omega)}^p dr
%\leq C(\Omega, n,p, s) \|f\|_{W^{s,p}(\Omega)}^p.
%$$
We can find $U\subset\subset \tilde{U}\subset\subset\Omega$ so that $\tilde{U}$ has smooth boundary
As $f_k\rightharpoonup f$ we obtain that $f_k$ forms a bounded sequence in $W^{s,n/s}(\tilde{U})$. We denote, $H_r:=\{x\in \rn\ |\ \pi_i(x)=r\}$ where $\pi_i$ is the projection to $x_i$ axis. 
With the help of Fatou's Lemma and Lemma~\ref{lm-2}, 
$$
\int_{\mathbb{R}} \liminf_{k\to\infty} \|f_k\|_{W^{s,n/s}(H_r\cap \tilde{U})}^{n/s} dr\leq \liminf_{k\to\infty} \int_{\mathbb{R}} \|f_k\|_{W^{s,n/s}(H_r\cap \tilde{U})}^{n/s} dr 
\leq C \liminf_{k\to\infty} \|f_k\|_{W^{s,n/s}(\tilde{U})}<\infty.
$$
It follows that, for $\mathcal{H}^1$-a.e. $r$, we have that $\liminf \|f_k\|_{W^{s,n/s}(U\cap H_r)}<\infty$. Thus for each such $H_r$ there is a subsequence (for each hyperplane it could be different subsequence) 
$f_{k_m}\rightharpoonup f$ in $W^{s,n/s}(H_r\cap U)$. As $s\tfrac{n}{s}>n-1$ it follows that $f_{k_m}$ converges weakly to  $f$ in some $C^{0,\alpha}(H_r\cap U)$ and thus uniformly on $H_r\cap U$ with the help of Lemma \ref{cpt-embedding-holder}.

Let $\epsilon>0$. 
%By Lemma \ref{imagepositive} we know that $f(U)$ is the bounded set surrounded by $f(\partial U)$. 
Using the observation of uniform convergence on a.e. hyperplane for hyperplanes parallel to coordinate axes and a simple selection argument we can obtain a domain $U_I\subset U$ 
%with $|U\setminus U_I|<\epsilon$ 
so that boundary of $U_I$ is a subset of finitely many hyperplanes parallel to coordinate axes and there is a subsequence $f_{k_j}$ which converges uniformly to $f$ on $\partial U_I$. 
Moreover, by choosing $U_I$ really close to $U$ we can assume that 
$$
|f(U)\setminus f(U_I)|<\epsilon.
$$ 
Recall that $|f(\partial U)|=0$ and let $\delta>0$ be small enough so that 
$$
|f(\partial U)+B(0,\delta)|<\epsilon.
$$ 
 Using uniform continuity of $f$ on $\overline{U}$ we can assume that $U_I$ is so close $U$ so that 
 $f(\partial U_I)\subset f(\partial U)+B(0,\delta)$ and hence 
$$
|f(\partial U_I)|\leq |f(\partial U)+B(0,\delta)|<\epsilon. 
$$
As $f_{k_j}$ converges uniformly to $f$ on $\partial U_I$ we can assume that $\|f-f_{k_j}\|_{\infty}<\delta$ on $\partial U_I$ for $j$ big enough and we claim that this implies 
$$
A_{\delta}:=\bigl\{y\in f(U_I):\ \dist(y,\partial f(U_I))>\delta\bigr\}\subset f_{k_j}(U_I).
$$
Indeed, from $y\in f(U_I)$ we obtain $\deg(f,U_I,y)>0$ as $f$ is sense-preserving and using Theorem \ref{aprox} (recall $\|f-f_{k_j}\|_{L^\infty(U_I)}<\delta<\dist(y,\partial f(U_I))$) we obtain $\deg(f_{k_j},U_I,y)=\deg(f,U_I,y)>0$ and hence $y\in f_{k_j}(U_I)$.
It is not difficult to see that by choosing $\delta$ small enough we can achieve that 
$$
|f(U_I)\setminus f_{k_j}(U_I)|\leq |f(U_I)\setminus A_{\delta}|
\leq |f(\partial U_I)|+|\bigl\{x\in f(U_I):\ 0<\dist(x,\partial f(U_I))\leq\delta\bigr\}| 
<2\epsilon. 
$$
It follows that 
$$
\liminf_{j} |f_{k_j}(U)|
\geq \liminf_{j} |f_{k_j}(U_I)|
\geq |f(U)|-3\epsilon. 
$$

To get the opposite inequality, we approximate $U$ from outside by $U_O\supset U$, we follow the same steps as before and using moreover $|f(\partial U)|=0$ we get 
$$
\limsup_j |f_{k_j}(U)|\leq \limsup_j |f_{k_j}(U_O)|\leq |f(U)|+3\epsilon. 
$$
As from each subsequence of $f_k$ we can select such subsequence $f_{k_j}$ we obtain \eqref{haha}.
\end{proof}

\prt{Theorem}
\begin{proclaim}\label{th-CN-stability}
Let $\Omega\subset \rn$ be a bounded domain, $s\in \left(\frac{n-1}{n},1\right)$, and $p\in \left[\frac{n}{s},\infty\right)$. 
Assume that a sequence of mappings $\{f_k\}_k$, converges to $f$ weakly in $W^{s,p}(\Omega,\rn)$, assume that $\J_{f}>0$ and that $ \J_{f_k}>0$ for all $k$. %{\color{blue} %Additionally, in case of $sp=n$, 
Let $U\subset\subset \Omega$ be an open set with $|\partial f(U)|=0$. Assume that all $f_k$ do satisfy the Ciarlet-Ne\v{c}as condition
$$
\J_{f_k}(U)=|f_k(U)|, \text{ then }\J_{f}(U)=|f(U)|. 
$$
\end{proclaim}
\begin{proof}
Let us choose an $\varepsilon>0$. Then there is some $\psi\in C_c^\infty (\rn)$ with $\psi \leq \chi_U$ such that
\begin{equation}\label{eq2}
\J_f(\chi_U)
\leq \J_f(\psi) +\varepsilon.
\end{equation}
Since $sp\geq n$, we get by \cite[Theorem~1.8, Theorem~1.11]{DiVa} that $W^{s,p}(\Omega)$ is compactly embedded in $W^{\frac{n-1}{n},n}(\Omega)$. So, weak convergence of $f_k$ in $W^{s,p}(\Omega)$ implies, up to a subsequence, that $f_k\to f$ in $W^{\frac{n-1}{n},n}(\Omega)$. From Theorem~\ref{BN-theorem}, we then have
\begin{equation}\label{eq1}
 \J_{f_k}(\psi) \to  \J_f(\psi) \quad \text{for all } \psi\in C_c^\infty (\rn).
\end{equation}
We combine \eqref{eq1} and \eqref{eq2} to conclude that
$$
\J_f(\chi_U)
\leq \lim_{k\to\infty} \J_{f_k}(\psi) +\varepsilon.
$$
The definition of $\J_{f_k}(\chi_U)$ implies $\J_{f_k}(\psi) \leq \J_{f_k}(\chi_U)$. So, we have
$$
\J_f(\chi_U)
\leq \lim_{k\to\infty} \J_{f_k}(\chi_U) +\varepsilon.
$$
To prove the opposite inequality we use the fact that $U\subset\subset\Omega$ and we choose $\phi\in C_c(\Omega)$ so that $\chi_{U}\leq \phi$ and analogously as above 
$$
\J_f(\chi_U)\geq \J_f(\phi)-\epsilon=\lim_{k\to\infty}\J_{f_k}(\phi)-\epsilon\geq \lim_{k\to\infty} \J_{f_k}(\chi_U) -\varepsilon.
$$
It follows that
$$
\lim_{k\to\infty} \J_{f_k}(\chi_U)=\J_f(\chi_U).
$$ 
By Theorem \ref{sense} we know that assumption of Lemma \ref{lm-1} are satisfied and thus Lemma \ref{lm-1} gives us $\lim_{k\to\infty} |f_k(U)|=|f(U)|$ and our conclusion follows. 
\end{proof}

\section{Counterexample to Lusin $(N)$ condition in $W^{s,p}$}

In this section we prove Theorem~\ref{exampleJarda}. Let us fix $0<\alpha<1$. We construct an $\alpha$-H\"older continuous homeomorphism which maps a Cantor type set with zero measure $C_A$ onto a Cantor type set with nonzero measure $C_B$. 
%We have shown that the mapping $f$ is $\alpha$-H\"older continuous homeomorphism which maps $\mathcal{C}_A$ onto $\mathcal{C}_B$ and thus fails the Lusin (N) condition. 
It is easy to see that $C^{0,\alpha}(Q_0)\subset W^{s,p}(Q_0)$ whenever $\alpha>s$, consequently, for any $s\in(0,1)$ and $p\in[1,\infty)$, we can construct such homeomorphism.

 The construction is inspired by \cite[Section 4.3]{HK}, more details could be found there. Let $0<\alpha<1$ and set 
$$
A=2^{\frac{1-\alpha}{\alpha}},\text{ i.e. } 2=(2A)^{\alpha}. 
$$

\subsection{Construction of Cantor type sets}\label{defCantor} 
We first give two Cantor set constructions in $(-1,1)^n$. %Our mapping $f$ will
%be defined as a limit of a sequence of piecewise continuously differentiable homeomorphisms
%$f_k:(-1,1)^n\to(-1,1)^n$, where each $f_k$ maps the $k$-th step of the first Cantor
%set construction onto the second one. Then the limit mapping $f$ maps the first Cantor set onto the
%second one. 

By $\V$ we denote the set of $2^{n}$ vertices of the cube $[-1,1]^{n}$.
The sets $\V^k=\V\times\hdots\times \V$, $k\in\en$, will serve as the sets of
indices for our construction. 
Let us denote 
\eqn{defa}
%\[a_k=A^{-k}~\text{and}~b_k=\frac{1}{2}\left(1-\frac{1}{A}\right)+A^{-k}.\]
$$
a_k=A^{-k}\ \text{ and }\ b_k=\frac{1}{2}\Bigl(1-\frac{1}{A}\Bigr)+A^{-k}.
$$
Set $z_0=\z_0=0$ and let us define 
\eqn{defr}
$$r_k=a_k 2^{-k}\ \text{ and }\ \r_k=b_k 2^{-k} .$$
It follows that $(-1,1)^{n}=Q(z_0,r_0)$ and further we proceed by induction. 
For $\ve=[v_1,\hdots, v_k]\in \V^{k}$ we denote $\w=[v_1,\hdots,v_{k-1}]$ and we define 
$$z_{\ve}=z_{\w}+\frac{1}{2} r_{k-1} v_k=z_0+\frac{1}{2}\sum_{j=1}^k r_{j-1}v_j,$$
$$Q'_{\ve}=Q(z_{\ve},\tfrac{r_{k-1}}{2})\ \text{ and }\ Q_{\ve}=Q(z_{\ve},r_k).$$
%Formally we should write $\w(\ve)$ instead of $\w$ but for the simplification of the notation we will avoid this.

\unitlength=0.7mm
\begin{picture}(110,50)(-50,5)
%\circle{10}
\put(10,10){\framebox(40,40){}}
\put(30,10){\line(0,1){40}}
\put(10,30){\line(1,0){40}}
\put(12,12){\framebox(16,16){}}
%\put(20,32){\makebox(5,5){$R_{[-1,1]}$}}
\put(12,32){\framebox(16,16){}}
\put(32,12){\framebox(16,16){}}
\put(32,32){\framebox(16,16){}}
%end of first cubeanalogs of the computation above together with
\put(60,10){\framebox(40,40){}}
\put(80,10){\line(0,1){40}}
\put(60,30){\line(1,0){40}}
\put(62,12){\framebox(16,16){}}%1. rectangle
\put(70,12){\line(0,1){16}}
\put(62,20){\line(1,0){16}}
\put(63,13){\framebox(6,6){}}
\put(71,13){\framebox(6,6){}}
\put(63,21){\framebox(6,6){}}
\put(71,21){\framebox(6,6){}}
\put(62,32){\framebox(16,16){}}%2. rectangle
\put(70,32){\line(0,1){16}}
\put(62,40){\line(1,0){16}}
\put(63,33){\framebox(6,6){}}
\put(71,33){\framebox(6,6){}}
\put(63,41){\framebox(6,6){}}
\put(71,41){\framebox(6,6){}}
\put(82,12){\framebox(16,16){}}%3. rectangle
\put(90,12){\line(0,1){16}}
\put(82,20){\line(1,0){16}}
\put(83,13){\framebox(6,6){}}
\put(91,13){\framebox(6,6){}}
\put(83,21){\framebox(6,6){}}
\put(91,21){\framebox(6,6){}}
\put(82,32){\framebox(16,16){}}%4. rectangle
\put(90,32){\line(0,1){16}}
\put(82,40){\line(1,0){16}}
\put(83,33){\framebox(6,6){}}
\put(91,33){\framebox(6,6){}}
\put(83,41){\framebox(6,6){}}
\put(91,41){\framebox(6,6){}}
\end{picture}

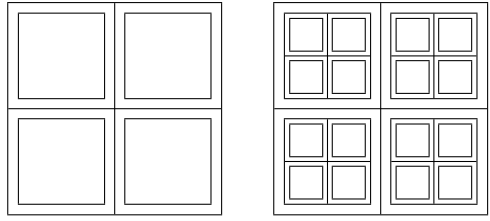
\centerline{{\bf Fig. 1.} Cubes $Q_{\ve}$ and $Q'_{\ve}$ for $\ve\in \V^1$ and $\ve\in \V^2$.}
\vskip 10pt

The number of the cubes $\{Q_{\ve}:\ \ve\in\V^k\}$ is $2^{nk}$. It is not difficult to find out that the resulting Cantor set
$$\bigcap_{k=1}^{\infty} \bigcup_{\ve\in\V^k} Q_{\ve}=:C_A=C_a\times \hdots \times C_a$$
is a product of $n$ Cantor sets in $\er$. Moreover, $\mira(C_A)=0$ since 
$$\mira\bigl(\bigcup_{\ve\in\V^k} Q_{\ve}\bigr)=2^{nk} (2 a_k 2^{-k})^{n}\overset{k\to\infty}{\to} 0.$$

Analogously we define 
$$\z_{\ve}=\z_{\w}+\frac{1}{2} \r_{k-1} v_k=\z_0+\frac{1}{2}\sum_{j=1}^k \r_{j-1}v_j,$$
$$\Q'_{\ve}=Q(\z_{\ve},\tfrac{\r_{k-1}}{2})\text{ and }\Q_{\ve}=Q(\z_{\ve},\r_k).$$
The resulting Cantor set 
$$\bigcap_{k=1}^{\infty} \bigcup_{\ve\in\V^k} \Q_{\ve}=:C_B=C_b\times \hdots \times C_b$$
satisfies $\mira(C_B)>0$ since $\lim_{k\to\infty}b_k>0$. 
It remains to find a homeomorphism $f$ which maps $C_A$ onto $C_B$ and satisfies our assumptions. 
From $\mira(C_A)=0$ and $\mira(C_B)>0$ we will obtain that $f$ does not satisfy the $(N)$ condition. 

\unitlength=0.6mm
\begin{picture}(115,60)(-75,5)
\put(10,20){\framebox(30,30){}}%domain
\put(20,30){\framebox(10,10){}}
%\put(15,15){\dashbox(20,40){}}
\put(10,20){\line(1,1){10}}
\put(40,20){\line(-1,1){10}}
\put(10,50){\line(1,-1){10}}
\put(40,50){\line(-1,-1){10}}
\put(32,53){\makebox(5,5){$Q'_{\ve}$}}
\put(23,42){\makebox(5,5){$Q_{\ve}$}}
%\put(13,32){\makebox(5,5){$A$}}
%\put(33,32){\makebox(5,5){$A$}}
%\put(23,12){\makebox(5,5){$B$}}
%\put(23,52){\makebox(5,5){$B$}}
\put(42,35){\vector(1,0){10}}%sipka
\put(45,36){\makebox(5,5){$f_k$}}
\put(55,15){\framebox(40,40){}}%target
\put(65,25){\framebox(20,20){}}
%\put(65,25){\dashbox(50,20){}}
\put(55,15){\line(1,1){10}}
\put(95,15){\line(-1,1){10}}
\put(55,55){\line(1,-1){10}}
\put(95,55){\line(-1,-1){10}}
\put(80,58){\makebox(5,5){$\Q'_{\ve}$}}
\put(75,48){\makebox(5,5){$\Q_{\ve}$}}
\end{picture}

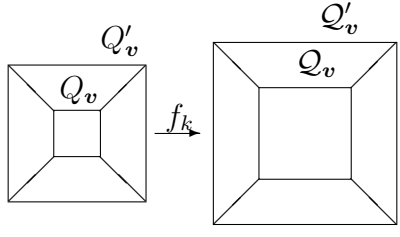
\centerline{{\bf Fig. 2.} The transformation of $Q'_{\ve}\setminus Q_{\ve}$ onto $\Q'_{\ve}\setminus \Q_{\ve}$}\label{pic}
\vskip 10pt

\subsection{Construction of a mapping $f$}\label{defmap} 

%Again we will proceed by induction and 
We will find a sequence of homeomorphisms 
$f_k:(-1,1)^{n}\to(-1,1)^{n}$. We set $f_0(x)=x$ and we proceed by induction. 
We will give a 
mapping $f_1$ which stretches each cube $Q_{v}$, $v\in\V^1$ , homogeneously so that $f_1(Q_{v})$
equals $\Q_{v}$. On the annulus $P_v\setminus Q_v$, $f_1$ is defined to be an appropriate radial map
with respect to $z_v$ and $\z_v$ in the image in order to make $f_1$ a homeomorphism. The
general step is the following: If $k > 1$, $f_k$ is defined as $f_{k-1}$ outside the union of all
cubes $Q_{\w}$, $\w\in \V^{k-1}$. Further, $f_k$ remains equal to $f_{k-1}$ at the centers of cubes
$Q_{\ve}$, $\ve\in \V^k$. Then $f_k$ stretches each cube $Q_{\ve}$, $\ve\in \V^k$, homogeneously so that
$f_k(Q_{\ve})$ equals $\Q_{\ve}$. On the annulus $P_{\ve}\setminus Q_{\ve}$, $f_k$ is defined to be an appropriate radial
map with respect to $z_{\ve}$ in preimage and $\z_{\ve}$ in image to make $f_k$ a homeomorphism 
(see Fig. 2). %Notice that the Jacobian determinant $J_{f_k}(x)$ will be strictly
%positive almost everywhere in $(-1,1)^n$.

In this proof we use the notation $\|x\|$ for the supremum norm of $x\in\rn$. 
The mappings $f_k$, $k\in\en$, are formally defined as 
\eqn{defg}
$$
f_k(x)=
\begin{cases}
f_{k-1}(x)&\text{for }x\notin\bigcup_{\ve\in\V^k}Q'_{\ve}\\
f_{k-1}(z_{\ve})+(\alpha_k\|x-z_{\ve}\|+\beta_k)\frac{x-z_{\ve}}{\|x-z_{\ve}\|}
&\text{for }x\in Q'_{\ve}\setminus Q_{\ve},\ \ve\in\V^k\\
f_{k-1}(z_{\ve})+\frac{\r_k}{r_k}(x-z_{\ve})&\text{for }x\in Q_{\ve},\ \ve\in\V^k\\
\end{cases}
$$
where the constants $\alpha_k$ and $\beta_k$ are given by 
\eqn{defalpha}
$$
\alpha_k r_k+\beta_k=\r_k\text{ and }\alpha_k \tfrac{r_{k-1}}{2}+\beta_k=\tfrac{\r_{k-1}}{2} .
$$
It is not difficult to find out that each $f_k$ is a homeomorphism and maps 
$$\bigcup_{\ve\in\V^k} Q_{\ve} \text{ onto }\bigcup_{\ve\in\V^k} \Q_{\ve} .$$
The limit $f(x)=\lim_{k\to\infty}f_k(x)$ is clearly one-to-one and continuous and therefore a homeomorphism. 
Moreover, it is easy to see that $f$ is differentiable almost everywhere %, absolutely continuous on almost all lines parallel to coordinate axes 
and maps $C_A$ onto $C_B$. %Moreover, the pointwise Jacobian is positive a.e. and hence our map is a mapping of finite distortion if the derivative is integrable. 

\subsection{H\"older continuity of $f$} 

Let $k\in\en$ and $\ve\in \V^k$. 
We need to estimate $Df_k(x)$ %and $J_g(x)$ 
in the interior of the annulus $Q'_{\ve}\setminus Q_{\ve}$. Since 
$$
f_k(x)=f_k(z_{\ve})+(\alpha_k \|x-z_{\ve}\|+\beta_k)\frac{x-z_{\ve}}{\|x-z_{\ve}\|}
$$
there, we can almost apply Lemma \ref{radial} to compute its derivative. 
The difference is that instead of the center $0$ we use $z_{\ve}$ and also here we have supremum norm while in the statement of Lemma \ref{radial} we have used an euclidean norm. But these norms are equivalent up to a bilipschitz change of variables and therefore we can compute the norm of the derivative and the result is comparable up to a multiplicative constant $C(n)$. 
From Lemma \ref{radial}, $r_k\sim r_{k-1}$, $\r_k\sim \r_{k-1}$ \eqref{defalpha}, \eqref{defr} and \eqref{defa} 
we now obtain
$$
|Df_k(x)|\sim \max\Bigl\{\frac{\r_k}{r_k}, \alpha_k\Bigr\}\sim 
\max\Bigl\{\frac{b_k}{a_k}, \frac{b_{k-1}-b_k}{a_{k-1}-a_k}\Bigr\}\leq C A^k.
$$
Since $f_k=f_{k-1}$ outside of $Q'_{\ve}$ it is easy to show by induction that each $f_k$ is globally $C A^k$ lipschitz continuous. 
%Moreover, we can estimate 
%$$\mira(Q'_{\ve}\setminus Q_{\ve})= (r_{k-1})^{n}-(2 r_k)^{n}\sim 
%2^{-kn}\Bigl(\frac{1}{(k-1)^n}-\frac{1}{k^n}\Bigr)\sim 2^{-kn}\frac{1}{k^{n+1}}$$
%and we have $2^{kn}$ annuli like that. Therefore, as $p<n$ 

Let $x,y\in(0,1)^n$. We find $k\in\en$ such that 
\begin{equation}\label{dist}	
    \sqrt{n}r_{k+1}\leq\|x-y\|\leq \sqrt{n}r_{k}.
\end{equation}
Set 
\[t_1=\sup\Bigl\{t\leq0:tf(x)-(1-t)f(y)\notin \bigcup_{\ve\in\V^{k+1}} \Q_{\ve}
%\bigcup_{\tilde{Q}\in\tilde{\mathcal{Q}}^\prime_{k+1}}\tilde{Q}
\Bigr\}
\]
and
\[t_2=\inf\Bigl\{t\geq 1:tf(x)-(1-t)f(y)\notin \bigcup_{\ve\in\V^{k+1}} \Q_{\ve}
\Bigr\}.\]
We also define
\[p_1=t_1f(x)-(1-t_1)f(y)~\text{and}~p_2=t_2f(x)-(1-t_2)f(y).\]
The points $p_1$ and $p_2$ are on the line connecting $f(x)$ and $f(y)$ and $\|p_1-p_2\|\geq \|f(x)-f(y)\|$. We find $x_0$ and $y_0$ such that $f(x_0)=p_1$ and $f(y_0)=p_2$. 
Note that either $x_0\in\partial(\bigcup_{\ve\in\V^{k+1}} Q_{\ve}  )$ or $x_0=x$ and similar thing holds for $y_0$, hence in both cases $f(x_0)=f_k(x_0)$ and $f(y_0)=f_k(y_0)$.

We have that
%\begin{align*}
$$
\|f(x)-f(y)\|
\leq\|f(x_0)-f(y_0)\|
=\|f_k(x_0)-f_k(y_0)\|. 
$$
%\end{align*}
Since $f_k$ is $CA^k$-lipschitz continuous we have that
\begin{align*}
\|f_k(x_0)-f_k(y_0)\|
&\leq CA^k\|x_0-y_0\|\\
&\leq CA^k(\|x-x_0\|+\|y-y_0\|+\|x-y\|).
\end{align*}
From the definition of $x_0$ it follows that either $x=x_0$ or both $x$ and $x_0$ are in the same cube $Q_{\ve}$ for $\ve\in \V^{k+1}$,  
%$Q\in\tilde{\mathcal{Q}}_{k+1}$
hence $\|x-x_0\|\leq \sqrt{n}r_k/2$. The same holds for $y$ and $y_0$. This, together with (\ref{dist}) gives us
$$
CA^k(\|x-x_0\|+\|y-y_0\|+\|x-y\|)
\leq CA^k(\sqrt{n}r_k+\sqrt{n}r_k)
=2C\sqrt{n}A^k r_k
$$
It follows from $2=(2A)^{\alpha}$ and (\ref{defr}) that 
\begin{align*}
\|f(x)-f(y)\|\leq C  A^k r_k&=C 2^{-k-1}=C (2A)^{(-k-1)\log_{2A}2}\\
&=C r_{k+1}^\alpha\leq C \|x-y\|^\alpha.\\
\end{align*}
\qed

%2C\sqrt{n}A^kd_k
%&=4C\sqrt{n}2^{-k-1}\\
%&=4C\sqrt{n}(2A)^{(-k-1)\log_{2A}2}
%\end{align*}
%Again, we use (\ref{Ak}) together with the fact that  and get
%\begin{align*}
%4C\sqrt{n}(2A)^{(-k-1)\log_{2A}2}
%&=4C\sqrt{n}d_{k+1}^\alpha\\
%&\leq 4C\sqrt{n}\|x-y\|^\alpha.
%\end{align*}
%Hence, it follows that $\|f(x)-f(y)\|\leq C\|x-y\|^\alpha.$

\section{Example of mapping with $J_f>0$ a.e. which is not sense preserving in $W^{s,p}$}
In this section we prove Theorem~\ref{exampleKaushik} and we are using some ideas of counterexamples from \cite{KKM}. Again we fix $0<\alpha<1$ and it is enough to construct a mapping which is $\alpha$-H\"older continuous. We fix $A>0$ and $0<B<1$ small enough so that
\eqn{chooseAB}
$$
B+A<(1+A)(1-\alpha)\text{, i.e. }1-B\geq \alpha(1+A). 
$$
We set 
\eqn{defab}
$$
a_k=2^{-Ak}\ \text{ and }\ b_k=2^{Bk}.
$$
and as in previous section we define 
$$
r_k=a_k 2^{-k}\ \text{ and }\ \r_k=b_k 2^{-k} .
$$

\subsection{Definition of Cantor sets and mapping}
We define $C_A$ in the same way as in Section \ref{defCantor} using $a_k$ from \eqref{defab}. Analogously to Section \ref{defCantor} we also define $\Q_{\ve}$ and $\Q_{\ve}'$ using the value of $b_k$ from \eqref{defab}. The main difference here is that $b_k>1$ so $\Q_{\ve}$ is actually bigger than $\Q_{\ve}'$. However as $0<B<1$ and $\r_k=b_k 2^{-k}$ it is easy to see that 
$\r_k\to 0$ so we still have $\lim_{k\to\infty}\diam \Q_{\ve}=0$. 

%\vskip -10pt

\begin{tikzpicture}[line cap=round,line join=round,>=triangle 45,x=0.5cm,y=0.5cm]
\clip(-6.428,-13.198) rectangle (30.48,3.404);
\draw [line width=0.4pt,dash pattern=on 5pt off 5pt] (0.,0.)-- (4.,0.);
\draw [line width=0.4pt,dash pattern=on 5pt off 5pt] (4.,0.)-- (4.,-4.);
\draw [line width=0.4pt,dash pattern=on 5pt off 5pt] (4.,-4.)-- (0.,-4.);
\draw [line width=0.4pt,dash pattern=on 5pt off 5pt] (0.,-4.)-- (0.,0.);
\draw [line width=0.4pt] (-1.,1.)-- (-1.,-5.);
\draw [line width=0.4pt] (-1.,-5.)-- (5.,-5.);
\draw [line width=0.4pt] (5.,-5.)-- (5.,1.);
\draw [line width=0.4pt] (5.,1.)-- (-1.,1.);
\draw [->,line width=0.4pt] (6.42,-2.02) -- (9.28,-2.);
\draw [line width=0.4pt,dash pattern=on 5pt off 5pt] (11.,2.)-- (11.,-6.);
\draw [line width=0.4pt,dash pattern=on 5pt off 5pt] (11.,-6.)-- (19.,-6.);
\draw [line width=0.4pt,dash pattern=on 5pt off 5pt] (19.,-6.)-- (19.,2.);
\draw [line width=0.4pt,dash pattern=on 5pt off 5pt] (19.,2.)-- (11.,2.);
\draw [line width=0.4pt] (12.5,0.5)-- (12.5,-4.5);
\draw [line width=0.4pt] (12.5,-4.5)-- (17.5,-4.5);
\draw [line width=0.4pt] (17.5,-4.5)-- (17.5,0.5);
\draw [line width=0.4pt] (17.5,0.5)-- (12.5,0.5);
\draw [->,line width=0.4pt] (0.1,0.) -- (-1.,1.);
\draw [->,line width=0.4pt] (2.02,0.) -- (2.,1.);
\draw [->,line width=0.4pt] (4.,0.) -- (5.,1.);
\draw [->,line width=0.4pt] (4.,-2.06) -- (5.,-2.06);
\draw [->,line width=0.4pt] (4.,-4.) -- (5.,-5.);
\draw [->,line width=0.4pt] (2.12,-4.) -- (2.12,-5.);
\draw [->,line width=0.4pt] (0.,-4.) -- (-1.,-5.);
\draw [->,line width=0.4pt] (0.,-1.92) -- (-1.,-1.94);
\draw [->,line width=0.4pt] (11.,2.) -- (12.5,0.42);
\draw [->,line width=0.4pt] (15.1,2.) -- (15.08,0.5);
\draw [->,line width=0.4pt] (19.,2.) -- (17.5,0.5);
\draw [->,line width=0.4pt] (19.,-2.06) -- (17.5,-2.08);
\draw [->,line width=0.4pt] (19.,-6.) -- (17.5,-4.5);
\draw [->,line width=0.4pt] (15.12,-6.) -- (15.08,-4.5);
\draw [->,line width=0.4pt] (11.,-6.) -- (12.56,-4.5);
\draw [->,line width=0.4pt] (11.,-2.02) -- (12.5,-2.);
\draw (7.418,-0.594) node[anchor=north west] {$f_k$};
\draw (3.106,2.266) node[anchor=north west] {$Q_{\boldsymbol v}'$};
\draw (2.38,0.154) node[anchor=north west] {$Q_{\boldsymbol v}$};
\draw (17.01,3.102) node[anchor=north west] {${\mathcal Q}_{\boldsymbol v}$};
\draw (15.69,0.638) node[anchor=north west] {${\mathcal Q}_{\boldsymbol v}'$};
\end{tikzpicture}
\vskip -100pt

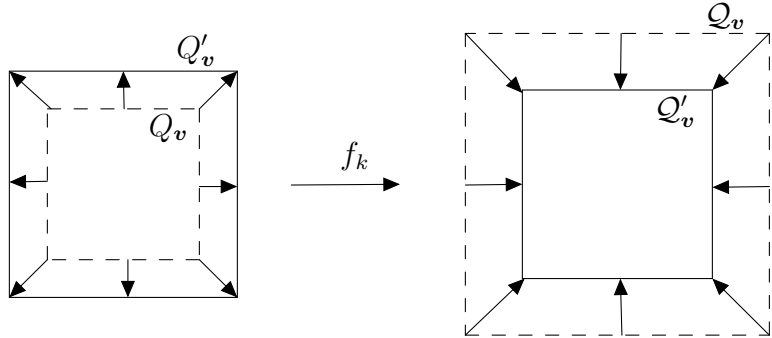
\centerline{{\bf Fig. 3.} Mapping $f_k$ maps $Q_{\ve}'\setminus Q_{\ve}$ onto $\Q_{\ve}\setminus \Q_{\ve}'$ and turns its orientation there.}
\vskip 10pt

We define our mapping $f$ as in Section \ref{defmap} using the same formula \eqref{defg} with $a_k$ and $b_k$ given by \eqref{defab}. 
The main difference here is that $Q_{\ve}\subset Q_{\ve}'$ but $\Q_{\ve}\supset \Q_{\ve}'$ so we are actually turning inside out (see Fig. 3) and it is easy to see that $J_{f_k}<0$ on $Q_{\ve}'\setminus Q_{\ve}$ and $J_{f_k}>0$ on $Q_{\ve}$. As $\lim_{k\to\infty}\diam \Q_{\ve}=0$ it is not difficult to see that the limiting mapping 
$f=\lim_{k\to\infty} f_k$ satisfies $f(x)=x$ on $\partial [-1,1]^n$ and it is in fact continuous (we even show H\"older continuity in the next section). 
Since $|C_A|=0$ it is easy to see that the classical derivative $\nabla f$ exists a.e. and at the end we have  $J_f=\det \nabla f<0$ a.e.

To obtain a mapping as promised in Theorem \ref{exampleKaushik} we define $\tilde{f}:Q_0\to\rn$ by
$$
\tilde{f}(x_1,x_2,\cdots,x_n)=f(-x_1,x_2,\cdots,x_n).
$$
This mapping clearly satisfies $J_{\tilde{f}}>0$ a.e. and $\tilde{f}(x_1,x_2,\cdots,x_n)=(-x_1,x_2,\cdots,x_n)$ for $x\in \partial Q_0$. 
Since $C^{0,\alpha}(Q_0)\subset W^{s,p}(Q_0)$ we obtain the required fractional Sobolev regularity once we prove that $f$ is $\alpha$-H\"older continuous.

\subsection{H\"older continuity of the map $f_k$}
Let $k\in\en$ and $\ve\in \V^k$. 
We need to estimate $Df_k(x)$ %and $J_g(x)$ 
in the interior of the annulus $Q'_{\ve}\setminus Q_{\ve}$. Since 
$$
f_k(x)=f_k(z_{\ve})+(\alpha_k \|x-z_{\ve}\|+\beta_k)\frac{x-z_{\ve}}{\|x-z_{\ve}\|}
$$
there, we can apply Lemma \ref{radial}, $r_k\sim r_{k-1}$, $\r_k\sim \r_{k-1}$, \eqref{defalpha} and \eqref{defab} 
and we obtain
$$
|Df_k(x)|\sim \max\Bigl\{\frac{\r_k}{r_k}, |\alpha_k|\Bigr\}\sim 
\max\Bigl\{\frac{b_k}{a_k}, \frac{|b_{k-1}-b_k|}{a_{k-1}-a_k}\Bigr\}\leq C 2^{(B+A)k}.
$$
Similarly for $x\in Q_{\ve}$ we have 
$$
|Df_k(x)|=\frac{\tilde{r}_k}{r_k}\leq C 2^{(B+A)k}.
$$
%Since $f_k=f_{k-1}$ outside of $Q_{\ve}'$ it is easy to see by induction that $f_k$ is globally $C 2^{(B+A)k}$ Lipschitz. 

For every $x,y\in Q'_{\ve}$ we have 
$$
|x-y|\leq \diam Q'_{\ve}\leq C r_k\leq C 2^{-k}2^{-Ak}
$$
and hence it is easy to see using \eqref{chooseAB} that 
\eqn{odhad}
$$
\begin{aligned}
|f_k(x)-f_k(y)|&\leq \sup_{ Q'_{\ve}} |Df_k|\cdot |x-y|\leq C 2^{(B+A)k}|x-y|\\
&\leq C 2^{(B+A)k} (2^{-k}2^{-Ak})^{1-\alpha} |x-y|^{\alpha}\leq C |x-y|^{\alpha}.\\
\end{aligned}
$$

\subsection{Global H\"older continuity of the map}  Let $x,y\in(-1,1)^n$ and find the smallest $k$ so that there is $Q'_{\w}\ni x, y$, $\w\in \V^{k-1}$. 
In the case $f(x)=f_k(x)$ and $f(x)=f_k(x)$ we know the required H\"older continuity from \eqref{odhad}. 

If this is not the case then without loss of generality $f(x)\neq f_k(x)$ which means that there is $\ve\in \V^{k}$ so that $x\in Q_{\ve}$. Since $k$ is the smallest possible $k$ we know that $y\notin Q'_{\ve}$ so it is not difficult to see that 
$$
|x-y|\geq \dist(Q_{\ve},\partial Q'_{\ve})\geq C 2^{-k}2^{-Ak}.
$$
We can now estimate using \eqref{odhad}
$$
\begin{aligned}
|f(x)-f(y)|&\leq |f_k(x)-f_k(y)|+|f_k(x)-f(x)|+|f_k(y)-f(y)|\\
&\leq C |x-y|^{\alpha}+ |f_k(x)-f(x)|+|f_k(y)-f(y)|.\\
\end{aligned}
$$
As $f=\lim f_k$ we can estimate using \eqref{defab}, $0<B<1$ and \eqref{chooseAB}
$$
\begin{aligned}
|f_k(x)-f(x)|&+|f_k(y)-f(y)|\leq 2\sum_{i=k}^{\infty}\|f_i-f_{i+1}\|_{\infty}
\leq C\sum_{i=k}^{\infty}\tilde{r_i}\\
&\leq C \sum_{i=k}^{\infty}2^{-i}2^{Bi}\leq C 2^{-k}2^{Bk}
\leq C (2^{-k}2^{-Ak} )^{\alpha}\leq C |x-y|^{\alpha}.\\
\end{aligned}
$$

\end{document}